\documentclass[12pt,a4paper]{article}
\usepackage[english]{babel}
\usepackage{amssymb}
\usepackage[T1]{fontenc}
\usepackage{amsmath}
\newtheorem{Th}{Theorem}
\newtheorem{Cor}{Corollary}
\newtheorem{Lem}{Lemma}
\newtheorem{RMQ}{Remark}

\newcommand{\F}{\mathcal{F}}
\newcommand{\R}{\mathbb{R}}
\newcommand{\Z}{\mathbb{Z}}

\newcommand{\N}{\mathbb{N}}
\newcommand{\M}{\mathcal{M}}

\newcommand{\G}{\mathcal{G}}

\renewcommand{\L}{\mathcal{L}}

\newcommand{\ds}{\displaystyle}
\newcounter{tictac}

\def\1{\,\rlap{\mbox{\small\rm 1}}\kern.15em 1}
\def\ind#1{\1_{#1}}
\def\build#1_#2^#3{\mathrel{\mathop{\kern 0pt#1}\limits_{#2}^{#3}}}
\def\tend#1#2{\build\hbox to 12mm{\rightarrowfill}_{#1\rightarrow #2}^{a.s.}}

\def\converge#1#2#3#4{\build\hbox to
#1mm{\rightarrowfill}_{#2\rightarrow #3}^{\hbox{\scriptsize #4}}}
\begin{document}
\title{Exact convergence rates in the central limit theorem for a class of martingales}
\author{M. El Machkouri and L. Ouchti}
\maketitle
\begin{abstract}
We give optimal convergence rates in the central limit theorem for
a large class of martingale difference sequences with bounded
third moments. The rates depend on the behaviour of the
conditional variances and for stationary sequences the rate
$n^{-1/2}\log n$ is reached. We give interesting examples of
martingales with unbounded increments which belong to the
considered class.

\vspace{9cm}\hspace{-0.7cm}
\end{abstract}
\thispagestyle{empty}
\newpage
\section{Introduction and notations}
The optimal rate of convergence in the central limit theorem for
independent random variables $(X_{i})_{i\in\Z}$ is well known to
be of order $n^{-1/2}$ as soon as the $X_{i}$'s are centered and
have uniformly bounded third moments (cf. Berry \cite{Berry} and
Esseen \cite{Esseen}). For dependent random variables the rate of
convergence was also most fully investigated but in many results
the rate is not better than $n^{-1/4}$. For example, Philipp
\cite{Philipp69} obtains a rate of $n^{-1/4}(\log n)^{3}$ for
uniformly mixing sequences, Landers and Rogge
\cite{Landers-Rogge76} obtain a rate of $n^{-1/4}(\log n)^{1/4}$
for a class of Markov chains and Sunklodas \cite{Sunklodas77}
obtains a rate of $n^{-1/4}\log n$ for strong mixing sequences.
However, Rio \cite{Rio96} has shown that the rate $n^{-1/2}$ is
reached for uniformly mixing sequences of bounded random variables
as soon as the sequence $(\phi_{p})_{p>0}$ of uniform mixing
coefficients satisfies $\sum_{p>0}p\phi_{p}<\infty$. Jan
\cite{Jan} established also a $n^{-1/2}$ rate of convergence in
the central limit theorem for bounded processes taking values in
$\R^{d}$ under some mixing conditions and recently, using a
modification of the proof in Rio \cite{Rio96}, Le Borgne and
P\`ene \cite{leborgne-Pene} obtained the rate $n^{-1/2}$ for
stationary processes satisfying a strong decorrelation hypothesis.
For bounded martingale difference sequences, Ibragimov
\cite{Ibrag63} has obtained the rate of $n^{-1/4}$ for some
stopping partial sums and Ouchti \cite{Ouchti} has extended
Ibragimov's result to a class of martingales which is related to
the one we are going to consider in this paper. Several results in
the rate of convergence for the martingale central limit theorem
have been obtained for the whole partial sums, one can refer to
Hall and Heyde \cite{Hall-Heyde} (section 3.6.), Chow and Teicher
\cite{Chow-Teicher} (Theorem 9.3.2), Kato \cite{Kato}, Bolthausen
\cite{Bolth2} and others. In fact, Kato \cite{Kato} obtains for
uniformly bounded variables the rate $n^{-1/2}(\log n)^{3}$ under
the strong assumption that the conditional variances are almost
surely constant. In this paper, we are most interested in results
by Bolthausen \cite{Bolth2} who obtained the better (in fact
optimal) rate $n^{-1/2}\log n$ under somewhat weakened conditions.
In this paper, we shall not pursue to improve the rate
$n^{-1/2}\log n$ but rather introduce a large class of martingales
which leads to it. Finally, note that El Machkouri and Voln\'y
\cite{EM-Volny-tll} have shown that the rate
of convergence in the central limit theorem can be arbitrary slow 
for stationary sequences of bounded (strong mixing) martingale 
difference random variables.\\
Let $n\geq 1$ be a fixed integer. We consider a finite sequence
$X=(X_{1},...,X_{n})$ of martingale difference random variables (
i.e. $X_{k}$ is $\F_{k}$-measurable and $E(X_{k}\vert\F_{k-1})=0$
a.s. where $(\F_{k})_{0\leq k\leq n}$ is an increasing filtration
and $\F_{0}$ is the trivial $\sigma$-algebra). In the sequel, we
are going to use the following notations
$$
\sigma_{k}^{2}(X)=E(X_{k}^{2}\vert\F_{k-1}),\quad
\tau_{k}^{2}(X)=E(X_{k}^{2}),\quad 1\leq k\leq n,
$$
$$
v_{n}^{2}(X)=\sum_{k=1}^{n}\tau_{k}^{2}(X)\quad\textrm{and}\quad
V_{n}^{2}(X)=\frac{1}{v_{n}^{2}(X)}\sum_{k=1}^{n}\sigma_{k}^{2}(X).
$$
We denote also $S_{n}(X)=X_{1}+X_{2}+...+X_{n}$. 
The central limit theorem established by Brown \cite{Brown71} 
and Dvoretzky \cite{Dvoretzky70} states that under some
Lindeberg type condition 
$$
\Delta_{n}(X)=\sup_{t\in\R}\big\vert\mu\left(S_{n}(X)/v_{n}(X)\leq
t\right)-\Phi(t)\big\vert\converge{10}{n}{+\infty}{}0.
$$
For more about central limit theorems for martingale difference 
sequences one can refer to Hall and Heyde \cite{Hall-Heyde}.  The rate of convergence of $\Delta_{n}(X)$ to zero was most fully investigated. Here, we focus on the following result by Bolthausen
\cite{Bolth2}.\\
\vspace{-0.2cm}
\\
\textbf{Theorem (Bolthausen, 82)}\label{TheobyBolthausen} {\em Let
$\gamma>0$ be fixed. There exists a constant $L(\gamma)>0$
depending only on $\gamma$ such that for all finite martingale
difference sequence $X=(X_{1},...,X_{n})$ satisfying
$V_{n}^{2}(X)=1$ a.s. and $\|X_{i}\|_{\infty}\leq\gamma$ then}
$$
\Delta_{n}(X)\leq L(\gamma)\left(\frac{n\log n}{v_{n}^{3}}\right).
$$
We are going to show that the method used by Bolthausen
\cite{Bolth2} in the proof of the theorem above can be extended to
a large class of unbounded martingale difference sequences. Note
that Bolthausen has already given extensions to unbounded
martingale difference sequences which conditional variances become
asymptotically nonrandom (cf. \cite{Bolth2}, Theorems 3 and 4) but
his assumptions cannot be compared directly with ours (cf.
condition ($\ref{class}$) below). So the results are
complementary.
\section{Main Results}
We introduce the following class of martingale difference
sequences: a sequence $X=(X_{1},...,X_{n})$ belongs to the class
$\M_{n}(\gamma)$ if $X$ is a martingale difference sequence with
respect to some increasing filtration $(\F_{k})_{0\leq k\leq n}$
such that for any $1\leq k\leq n$,
\begin{equation}\label{class}
E(\vert X_{k}\vert^{3}\vert\F_{k-1})\leq \gamma_{k}\,
E(X_{k}^{2}\vert\F_{k-1})\quad\textrm{a.s.}
\end{equation}
where $\gamma=(\gamma_{k})_{k}$ is a sequence of real numbers.\\
\\
Our first result is the following.
\begin{Th}\label{maintheorem}
Let $\gamma$ be a sequence of real numbers. There exists a
constant $L(\gamma)>0$ (not depending on $n$) such that for all
finite martingale difference sequence $X=(X_{1},...,X_{n})$ which
belongs to the class $\M_{n}(\gamma)$ and which satisfies
$V_{n}^{2}(X)=1$ a.s. then
$$
\Delta_{n}(X)\leq L(\gamma)\left(\frac{u_{n}\log n}{v_{n}}\right)
$$
where
$u_{n}=\vee_{k=1}^{n}\gamma_{k}$.
\end{Th}
\begin{RMQ} {\em Note that if the martingale difference sequence $X$ is
stationary then we obtain the rate $n^{-1/2}\log n$ which is
optimal (cf. Bolthausen \cite{Bolth2}).}
\end{RMQ}
As in Bolthausen \cite{Bolth2}, we derive the following result
when the restrictive condition $V_{n}^{2}(X)=1$ a.s. is relaxed.
\begin{Th}\label{cormaintheorem}
Let $\gamma$ be a sequence of real numbers. There exists a
constant $L(\gamma)>0$ (not depending on $n$) such that for all
finite martingale difference sequence $X=(X_{1},...,X_{n})$ which
belongs to the class $\M_{n}(\gamma)$ then
$$
\Delta_{n}(X)\leq L(\gamma)\left(\frac{u_{n}\log n}{v_{n}}
+\|V_{n}^{2}(X)-1\|^{1/2}_{\infty}\wedge\|V_{n}^{2}(X)-1\|^{1/3}_{1}\right)
$$
where $u_{n}$ is defined in Theorem $\ref{maintheorem}$.
\end{Th}
\textbf{Examples}. Let $X=(X_{1},..,X_{n})$ be a sequence of
martingale difference random variables such that $\sup_{1\leq
i\leq n}\|X_{i}\|_{\infty}\leq M<\infty$ and consider an arbitrary
sequence of independent random variables
$(\varepsilon_{1},...,\varepsilon_{n})$ with zero mean and finite
third moments which are also independent of $X$. One can notice
that the sequence $Y=(Y_{1},...,Y_{n})$ defined by
$Y_{k}=X_{k}+\varepsilon_{k}$ (resp. $Y_{k}=X_{k}\varepsilon_{k}$)
belongs to the class $\M_{n}(\gamma)$ where
$$
\gamma_{k}=4\left(M\vee\frac{E\vert\varepsilon_{k}\vert^{3}}{E\vert\varepsilon_{k}\vert^{2}}\right)
\quad\left(\textrm{resp.$\gamma_{k}=M\times\frac{E\vert\varepsilon_{k}\vert^{3}}{E\vert\varepsilon_{k}\vert^{2}}$}\right)
$$
Moreover if $V_{n}^{2}(X)=1$ a.s. (resp. $V_{n}^{2}(X)=1$ a.s. and
$(\varepsilon_{k})_{k}$ is stationary) then $V_{n}^{2}(Y)=1$ a.s.\\
\\
In the sequel, we are going to use the following lemma by
Bolthausen \cite{Bolth2}.
\begin{Lem}[Bolthausen, 82]\label{lemma-totalvariation}
Let $k\geq 0$ and $f:\R\to\R$ be a function which has $k$
derivatives $f^{(1)},...,f^{(k)}$ which together with $f$ belong
to $L^{1}(\mu)$. Assume that $f^{(k)}$ is of bounded variation
$\|f^{(k)}\|_{V}$, if $X$ is a random variable and if $\alpha_{1}\neq 0$
and $\alpha_{2}$ are two real numbers then
$$
\vert E
f^{(k)}(\alpha_{1}X+\alpha_{2})\vert\leq\|f^{(k)}\|_{V}\sup_{t\in\R}\vert\mu(X\leq
t)-\Phi(t)\vert+\vert\alpha_{1}\vert^{-(k+1)}\|f\|_{1}\sup_{x}\vert
\phi^{(k)}(x)\vert
$$
where $\phi(x)=(2\pi)^{-1/2}\exp(-x^{2}/2)$.
\end{Lem}
For any random variable $Z$ we denote
$\delta(Z)=\sup_{t\in\R}\vert\mu(Z\leq t)-\Phi(t)\vert$. We need
also the following extension of Lemma 1 in Bolthausen
\cite{Bolth2} which is of particular interest.
\begin{Lem}\label{lemma-inequalities}
Let $X$ and $Y$ be two real random variables. If there exist real
numbers $k>0$ and $r\geq 1$ such that $Y$ belongs to $L^{kr}(\mu)$
then there exist positive constants $c_{1}$ and $c_{2}$ such that
\begin{equation}\label{keyinequality1}
\delta(X+Y)\leq 2\delta(X)+c_{1}\,\|E\left(|Y|^{k}\vert
X\right)\|^{\frac{1}{k+1}}_{r}\wedge \|E\left(Y^{2}\vert
X\right)\|_{\infty}^{1/2}
\end{equation}
and
\begin{equation}\label{keyinequality2}
\delta(X)\leq 2\delta(X+Y)+c_{2}\,\|E\left(|Y|^{k}\vert
X\right)\|^{\frac{1}{k+1}}_{r}\wedge\| E\left(Y^{2}\vert
X\right)\|_{\infty}^{1/2}.
\end{equation}
\end{Lem}
The proofs of various central limit theorems for stationary
sequences of random variables are based on approximating the
partial sums of the process by martingales (see Gordin
\cite{Gordin69}, Voln\'y \cite{Volny93}). More precisely, if
$(f\circ T^{k})_{k}$ is a $p$-integrable stationary process where
$T:\Omega\to\Omega$ is a bijective, bimeasurable and
measure-preserving transformation (in fact, each stationary
process has such representation) then there exists necessary and
sufficient conditions (cf. Voln\'y \cite{Volny93}) for $f$ to be
equal to $m+g-g\circ T$ where $(m\circ T^{k})_{k}$ is a
$p$-integrable stationary martingale difference sequence and $g$
is a $p$-integrable function. In fact, such a decomposition can
hold also with $m$ and $g$ in some Orlicz space (see \cite{Volny-Orlicz}).
The term $g-g\circ T$ is called a coboundary. \\
The following theorem gives the rate of convergence in the central
limit theorem for stationary processes obtained from a martingale
difference sequence which is perturbed by a coboundary.
\begin{Th}\label{rate-coboundary}
Let $p>0$ be fixed and let $(f\circ T^{k})_{k}$ be a stationary
process. If there exist $m$ and $g$ in $L^{p}(\mu)$ such that
$(m\circ T^{k})_{k}$ is a martingale difference sequence and
$f=m+g-g\circ T$ then there exists a positive constant $c$ such
that
$$
\Delta_{n}(f)\leq
2\Delta_{n}(m)+\frac{2c\|g\|_{p}^{p/(p+1)}}{n^{p/2(p+1)}}.
$$
If $p=\infty$ then
$$
\Delta_{n}(f)\leq 2\Delta_{n}(m)+\frac{2c\|g\|_{\infty}}{n^{1/2}}.
$$
\end{Th}
Recently, Bosq \cite{Bosq2003} has shown that the condition
$\sum_{j\geq 0}j\vert\alpha_{j}\vert<\infty$ is sufficient to
obtain the optimal rate $n^{-1/2}$ for the linear process
$X_{k}=\sum_{j\geq 0}\alpha_{k-j}\varepsilon_{j}$ when
$(\varepsilon_{j})_{j}$ is a i.i.d. sequence with finite third
moment (Bosq established this result in the more general setting
of Hilbert spaces). We are going to give a convergence rate result
for linear processes when the $(\varepsilon_{j})_{j}$ are not
independent. By
\begin{equation}\label{def-linear-process}
X_{k}=\sum_{j\in\Z}\alpha_{k-j}\varepsilon_{j},\quad k\geq 1,
\end{equation}
we denote a stationary linear process, $(\varepsilon_{j})_{j}$ is
a stationary martingale difference sequence and $(\alpha_{j})_{j}$
are real numbers with $\sum_{j\in\Z}\alpha_{j}^{2}<\infty$.
\begin{Cor}\label{linear-process}
Let $\gamma$ be a sequence of real numbers. Consider the
linear process $(X_{k})_{k}$ defined by
$(\ref{def-linear-process})$ where $(\varepsilon_{j})_{j}$ is a
stationary martingale difference sequence such that
$E\vert\varepsilon_{0}\vert^{p}<\infty$ for some $p\geq 3$. Assume
that $\varepsilon=(\varepsilon_{1},...,\varepsilon_{n})$ belongs
to the class $\M_{n}(\gamma)$ and $V_{n}^{2}(\varepsilon)=1$ a.s. If
\begin{equation}\label{condition-decomposition}
\sum_{k=1}^{\infty}\left\{\bigg\vert\sum_{j\geq
k}\alpha_{j}\bigg\vert^{p}+\bigg\vert\sum_{j\leq
k}\alpha_{j}\bigg\vert^{p}\right\}<\infty
\end{equation}
then there exists a constant $L(\gamma)>0$ (not depending on $n$)
such that
$$
\Delta_{n}(X)\leq L(\gamma)\left(\frac{1}{n^{p/2(p+1)}}\right).
$$
If moreover $\varepsilon_{0}$ is a.s. bounded then
$$
\Delta_{n}(X)\leq L(\gamma)\left(\frac{\log n}{\sqrt{n}}\right).
$$
\end{Cor}
\begin{RMQ} {\em One can see that the condition $\sum_{j\in\Z}\vert
j\vert.\vert \alpha_{j}\vert<\infty$ is more restrictive than
($\ref{condition-decomposition}$).}
\end{RMQ}
\begin{RMQ} {\em Using Theorem $\ref{cormaintheorem}$, one can obtain a convergence rate result
for the linear process $(X_{k})_{k}$ when the condition
$V_{n}^{2}(\varepsilon)=1$ a.s. is relaxed.}
\end{RMQ}
\section{Proofs}
{\em Proof of Theorem $\ref{maintheorem}$}. Consider
$u=(u_{n})_{n}$ defined by $u_{n}=\vee_{k=1}^{n}\gamma_{k}$.
Clearly the class $\M_{n}(\gamma)$ is contained in the class
$\M_{n}(u)$. For any $(u,v)\in\R_{+}^{\N^{\ast}}\times\R_{+}^{\ast}$, we
consider the subclass
$$
\L_{n}(u,v)=\left\{X\in\M_{n}(u)\,\,\vert\,\,V_{n}^{2}(X)=1,\,\,v_{n}(X)=v\,\,\textrm{a.s.}
\right\}
$$
and we denote
$$
\beta_{n}(u,v)=\sup\left\{\Delta_{n}(X)\,\vert\,X\in\L_{n}(u,v)\right\}.
$$
In the sequel, we assume that $X=(X_{1},...,X_{n})$ belongs to
$\L_{n}(u,v)$, hence $X^{'}=(X_{1},...,X_{n-2},X_{n-1}+X_{n})$
belongs to $\L_{n-1}(4u,v)$ and consequently,
$$
\beta_{n}(u,v)\leq\beta_{n-1}(4u,v).
$$
Let $Z_{1},Z_{2},...,Z_{n}$ be independent identically distributed
standard normal variables independent of the $\sigma$-algebra
$\F_{n}$ (which contains the $\sigma$-algebra generated by
$X_{1},...,X_{n}$) and $\xi$ be an extra centered normal variable
with variance $\theta^{2}>1\vee 2u_{n}^{2}$ which is independent
of anything else. Noting that $\sum_{i=1}^{n}\sigma_{i}(X)Z_{i}/v$
is a standard normal random variable and according to Inequality
($\ref{keyinequality2}$) in Lemma \ref{lemma-inequalities},
\begin{equation}\label{inequalitybeginning}
\Delta_{n}(X)\leq
2\sup_{t\in\R}\big\vert\Gamma_{n}(t)\big\vert+\frac{3\theta}{v}.
\end{equation}
where
$$
\Gamma_{n}(t)\triangleq\mu\left(\left(S_{n}(X)+\xi\right)/v \leq
t\right)-\mu\left(\left(\sum_{i=1}^{n}\sigma_{i}(X)Z_{i}+\xi\right)/v\leq
t\right).
$$
For any integer $1\leq k\leq n$, we consider the following random
variables
$$
Y_{k}\triangleq\frac{1}{v}\sum_{i=1}^{k-1}X_{i},\quad
W_{k}\triangleq\frac{1}{v}\left(\sum_{i=k+1}^{n}\sigma_{i}(X)Z_{i}+\xi\right),
$$
$$
H_{k}\triangleq\frac{1}{v^{2}}\left(\sum_{i=k+1}^{n}\sigma_{i}^{2}(X)+\theta^{2}\right)\quad\textrm{and}\quad
T_{k}(t)\triangleq\frac{t-Y_{k}}{H_{k}},\,\,t\in\R
$$
with the usual convention
$\sum_{i=n+1}^{n}\sigma_{i}^{2}(X)=\sum_{i=n+1}^{n}\sigma_{i}(X)Z_{i}=0$
a.s. Moreover, one can notice that conditionned on
$\G_{k}=\sigma(X_{1},...,X_{n},Z_{k})$, the random variable
$W_{k}$ is centered normal with variance $H_{k}^{2}$. According to
the well known Lindeberg's decomposition (cf. \cite{Lindeberg}),
we have
\begin{align*}
\Gamma_{n}(t)
&=\sum_{k=1}^{n}\mu\left(Y_{k}+W_{k}+\frac{X_{k}}{v}\leq
t\right)-\mu\left(Y_{k}+W_{k}+\frac{\sigma_{k}(X)Z_{k}}{v}\leq
t\right)\\
&=\sum_{k=1}^{n}\mu\left(\frac{W_{k}}{H_{k}}\leq
T_{k}(t)-\frac{X_{k}}{vH_{k}}\right)-\mu\left(\frac{W_{k}}{H_{k}}\leq
T_{k}(t)-\frac{\sigma_{k}(X)Z_{k}}{vH_{k}}\right)\\
&=\sum_{k=1}^{n}E\left(E\left(\ind{\frac{W_{k}}{H_{k}}\leq
T_{k}(t)-\frac{X_{k}}{vH_{k}}}\,\vert\G_{k}\right)\right)-E\left(E\left(\ind{\frac{W_{k}}{H_{k}}\leq
T_{k}(t)-\frac{\sigma_{k}(X)Z_{k}}{vH_{k}}}\,\vert\G_{k}\right)\right)\\
&=\sum_{k=1}^{n}E\left(\Phi\left(T_{k}(t)-\frac{X_{k}}{vH_{k}}\right)\right)-E\left(\Phi\left(T_{k}(t)
-\frac{\sigma_{k}(X)Z_{k}}{vH_{k}}\right)\right)
\end{align*}
Now, for any integer $1\leq k\leq n$ and any random variable
$\zeta_{k}$, there exists a random variable
$\vert\varepsilon_{k}\vert<1$ a.s. such that
$$
\Phi\left(T_{k}(t)-\zeta_{k}\right)=\Phi\left(T_{k}(t)\right)-\zeta_{k}\Phi^{'}\left(T_{k}(t)\right)+\frac{\zeta_{k}^{2}}{2}\Phi^{''}\left(T_{k}(t)\right)
-\frac{\zeta_{k}^{3}}{6}\Phi^{'''}\left(T_{k}(t)-\varepsilon_{k}\zeta_{k}\right)\,\,\textrm{a.s.}
$$
So, we derive
\begin{align*}
\Gamma_{n}(t)=&\sum_{k=1}^{n}E\bigg\{\left(-\frac{X_{k}}{vH_{k}}+\frac{\sigma_{k}(X)Z_{k}}{vH_{k}}\right)\Phi^{'}(T_{k}(t))
+\left(\frac{X_{k}^{2}}{2v^{2}H_{k}^{2}}-\frac{\sigma_{k}^{2}(X)Z_{k}^{2}}{2v^{2}H_{k}^{2}}\right)\Phi^{''}(T_{k}(t))\\
&-\left(\frac{X_{k}^{3}}{6v^{3}H_{k}^{3}}\right)\Phi^{'''}\left(T_{k}(t)-\frac{\varepsilon_{k}X_{k}}{vH_{k}}\right)
+\left(\frac{\sigma_{k}^{3}(X)Z_{k}^{3}}{6v^{3}H_{k}^{3}}\right)\Phi^{'''}\left(T_{k}(t)-\frac{\varepsilon_{k}^{'}\sigma_{k}(X)Z_{k}}{vH_{k}}\right)\bigg\}.
\end{align*}
Since $V_{n}^{2}(X)=1$ a.s. we derive that $H_{k}$ and $T_{k}(t)$
are $\F_{k-1}$-measurable, hence
$$
\Gamma_{n}(t)
=\sum_{k=1}^{n}\frac{1}{6v^{3}}E\bigg\{-\frac{X_{k}^{3}}{H_{k}^{3}}\Phi^{'''}\left(T_{k}(t)-\frac{\varepsilon_{k}X_{k}}{vH_{k}}\right)
+\frac{\sigma_{k}^{3}(X)Z_{k}^{3}}{H_{k}^{3}}\Phi^{'''}\left(T_{k}(t)-\frac{\varepsilon_{k}^{'}\sigma_{k}(X)Z_{k}}{vH_{k}}\right)\bigg\}
$$
and consequently
\begin{equation}\label{inequalityinduction}
\big\vert\Gamma_{n}(t)\big\vert\leq \frac{1}{6v^{3}
}\left(S_{1}+S_{2}\right)
\end{equation}
where
$$
S_{1}=\sum_{k=1}^{n}E\bigg\{\frac{\vert
X_{k}\vert^{3}}{H_{k}^{3}}\bigg\vert\Phi^{'''}\left(T_{k}(t)-\frac{\varepsilon_{k}X_{k}}{vH_{k}}\right)\bigg\vert\bigg\}
$$
and
$$
S_{2}=\sum_{k=1}^{n}E\bigg\{\frac{\sigma_{k}^{3}(X)\vert
Z_{k}\vert^{3}}{H_{k}^{3}}\bigg\vert\Phi^{'''}\left(T_{k}(t)-\frac{\varepsilon_{k}^{'}\sigma_{k}(X)Z_{k}}{vH_{k}}\right)\bigg\vert\bigg\}.
$$
Consider the stopping times $\nu(j)_{j=0,..,n}$ defined by
$\nu(0)=0$, $\nu(n)=n$ and for any $1\leq j<n$
$$
\nu(j)=\inf\bigg\{k\geq
1\,\vert\,\sum_{i=1}^{k}\sigma_{i}^{2}(X)\geq\frac{jv^{2}}{n}\quad\textrm{a.s.}\bigg\}.
$$
Noting that
$\{1,...,n\}=\ds{\cup_{j=1}^{n}}\{\nu(j-1)+1,...,\nu(j)\}$ a.s. we
derive
$$
S_{1}=\sum_{j=1}^{n}E\bigg\{\sum_{k=\nu(j-1)+1}^{\nu(j)}\frac{\vert
X_{k}\vert^{3}}{H_{k}^{3}}\bigg\vert\Phi^{'''}\left(T_{k}(t)-\frac{\varepsilon_{k}X_{k}}{vH_{k}}\right)\bigg\vert\bigg\},
$$
moreover, for any $\nu(j-1)<k\leq\nu(j)$ we have
\begin{align*}
H_{k}^{2}&\geq\frac{1}{v^{2}}\left(\sum_{i=\nu(j)+1}^{n}\sigma_{i}^{2}(X)+\theta^{2}\right)\\
&=\frac{1}{v^{2}}\left(\sum_{i=1}^{n}\sigma_{i}^{2}(X)-\sum_{i=1}^{\nu(j)-1}\sigma_{i}^{2}(X)-\sigma_{\nu(j)}^{2}(X)+\theta^{2}\right)\\
&\geq\frac{1}{v^{2}}\left(v^{2}-\frac{jv^{2}}{n}-u_{n}^{2}+\theta^{2}\right)\\
&\triangleq m_{j}^{2}\quad\textrm{a.s.}
\end{align*}
Similarly,
\begin{align*}
H_{k}^{2}&\leq\frac{1}{v^{2}}\left(\sum_{i=\nu(j-1)+1}^{n}\sigma_{i}^{2}(X)+\theta^{2}\right)\\
&=\frac{1}{v^{2}}\left(\sum_{i=1}^{n}\sigma_{i}^{2}(X)-\sum_{i=1}^{\nu(j-1)}\sigma_{i}^{2}(X)+\theta^{2}\right)\\
&\leq\frac{1}{v^{2}}\left(v^{2}-\frac{(j-1)v^{2}}{n}+\theta^{2}\right)\\
&\triangleq M_{j}^{2}\quad\textrm{a.s.}
\end{align*}
Now, for any $\nu(j-1)<k\leq\nu(j)$ put
$$
R_{k}\triangleq\frac{1}{v}\sum_{i=\nu(j-1)+1}^{k-1}X_{i},\quad
A_{k}\triangleq\left\{\frac{\vert
R_{k}\vert}{m_{j}}\leq\frac{\vert
t-Y_{\nu(j-1)+1}\vert}{2M_{j}}\right\}
$$
and for any positive integer $q$ consider the real function
$\psi_{q}$ defined for any real $x$ by
$\psi_{q}(x)\triangleq\sup\{\vert\Phi^{'''}(y)\vert\,\,;\,\,y\geq\frac{\vert
x\vert}{2}-q\}$. In the other hand, on the set $A_{k}\cap\{\vert
X_{k}\vert\leq q\}$ we have
\begin{align*}
\big\vert T_{k}(t)-\frac{\varepsilon_{k}X_{k}}{vH_{k}}\big\vert
&=\big\vert\frac{t-Y_{\nu(j-1)+1}}{H_{k}}-\frac{R_{k}}{H_{k}}-\frac{\varepsilon_{k}X_{k}}{vH_{k}}\big\vert\\
&\geq\frac{\vert t-Y_{\nu(j-1)+1}\vert}{H_{k}}-\frac{\vert
R_{k}\vert}{H_{k}}-\frac{\vert X_{k}\vert}{vH_{k}}\\
&\geq\frac{\vert t-Y_{\nu(j-1)+1}\vert}{M_{j}}-\frac{\vert
R_{k}\vert}{m_{j}}-\frac{q}{\theta}\\
&\geq\frac{\vert
t-Y_{\nu(j-1)+1}\vert}{2M_{j}}-q\quad\textrm{a.s.}\quad\textrm{(since
$\theta\geq 1$).}
\end{align*}
Thus
$$
\bigg\vert\Phi^{'''}\left(T_{k}(t)-\frac{\varepsilon_{k}X_{k}}{vH_{k}}\right)\bigg\vert\ind{A_{k}\cap\vert
X_{k}\vert\leq
q}\leq\psi_{q}\left(\frac{t-Y_{\nu(j-1)+1}}{M_{j}}\right)\ind{A_{k}\cap\vert
X_{k}\vert\leq q}.
$$
So, for any $1\leq j\leq n$ we have
\begin{align*}
&E\bigg\{\sum_{k=\nu(j-1)+1}^{\nu(j)}\frac{|X_{k}|^{3}}{H_{k}^{3}}\bigg\vert
\Phi^{'''}\left(T_{k}(t)-\frac{\varepsilon_{k}X_{k}}{vH_{k}}\right)\bigg\vert\ind{A_{k}\cap\{\vert
X_{k}\vert\leq q\}}\bigg\}\\
&\leq E\bigg\{\sum_{k=\nu(j-1)+1}^{\nu(j)}\frac{\vert
X_{k}\vert^{3}}{H_{k}^{3}}\bigg\vert\psi_{q}\left(\frac{t-Y_{\nu(j-1)+1}}{M_{j}}\right)\bigg\vert\bigg\}\\
&=E\bigg\{E\bigg\{\sum_{k=\nu(j-1)+1}^{\nu(j)}\frac{\vert
X_{k}\vert^{3}}{H_{k}^{3}}\vert\F_{\nu(j-1)}\bigg\}\bigg\vert\psi_{q}\left(\frac{t-Y_{\nu(j-1)+1}}{M_{j}}\right)\bigg\vert\bigg\}\\
&=E\bigg\{E\bigg\{\sum_{k=\nu(j-1)+1}^{\nu(j)}E\left(\frac{\vert
X_{k}\vert^{3}}{H_{k}^{3}}\vert\F_{k-1}\right)\vert\F_{\nu(j-1)}\bigg\}\bigg\vert\psi_{q}\left(\frac{t-Y_{\nu(j-1)+1}}{M_{j}}\right)\bigg\vert\bigg\}\\
&\leq\frac{u_{n}}{m_{j}^{3}}
E\bigg\{E\bigg\{\sum_{k=\nu(j-1)+1}^{\nu(j)}\sigma_{k}^{2}(X)\vert\F_{\nu(j-1)}\bigg\}\bigg\vert\psi_{q}\left(\frac{t-Y_{\nu(j-1)+1}}{M_{j}}\right)\bigg\vert\bigg\}.
\end{align*}
Moreover, note that
\begin{align*}
\sum_{k=\nu(j-1)+1}^{\nu(j)}\sigma_{k}^{2}(X)&=\sum_{k=1}^{\nu(j)}\sigma_{k}^{2}(X)-
\sum_{k=1}^{\nu(j-1)}\sigma_{k}^{2}(X)\\
&\leq
\frac{(j+1)v^{2}}{n}-\frac{(j-1)v^{2}}{n}=\frac{2v^{2}}{n}\quad\textrm{a.s.}
\end{align*}
Using Lemma $\ref{lemma-totalvariation}$, noting that
$\|\psi_{q}\|_{\infty}\leq 1$ and keeping in mind the notation
$\delta(Z)\triangleq\sup_{t\in\R}\vert\mu(Z\leq t)-\Phi(t)\vert$
there exists a positive constant $c_{3}$ such that
$$
E\left\{\psi_{q}\left(\frac{t-Y_{\nu(j-1)+1}}{M_{j}}\right)\right\}
\leq\,\delta(Y_{\nu(j-1)+1})+c_{3}M_{j}.
$$
Now, using Lemma $\ref{lemma-inequalities}$ and the inequality
$$
E\left\{\left(\sum_{k=\nu(j-1)+1}^{n}X_{k}\right)^{2}\bigg\vert\F_{\nu(j-1)}\right\}\leq
v^{2}\left(1-\frac{j-1}{n}\right)\quad\textrm{a.s.}
$$
we obtain
\begin{align*}
\delta(Y_{\nu(j-1)+1})&\leq 2\,\delta(S_{n}(X)/v)
+c_{1}\bigg\| E\bigg\{\frac{1}{v^{2}}\left(\sum_{k=\nu(j-1)+1}^{n}X_{k}\right)^{2}\bigg\vert Y_{\nu(j-1)+1}\bigg\}\bigg\|_{\infty}^{1/2}\\
&=2\,\Delta_{n}(X)+c_{1}\bigg\| E\bigg\{\frac{1}{v^{2}}\left(\sum_{k=\nu(j-1)+1}^{n}X_{k}\right)^{2}\bigg\vert Y_{\nu(j-1)+1}\bigg\}\bigg\|_{\infty}^{1/2}\\
&\leq 2\,\beta_{n-1}(4u,v)+c_{1}\left(1-\frac{j-1}{n}\right)^{1/2}
\end{align*}
and so
$$
E\left\{\psi_{q}\left(\frac{t-Y_{\nu(j-1)+1}}{M_{j}}\right)\right\}\leq
2\,\beta_{n-1}(4u,v)+c_{1}\left(1-\frac{j-1}{n}\right)^{1/2}+c_{3}M_{j}.
$$
Using this estimate and the dominated convergence theorem, we
derive for any integer $1\leq j\leq n$,
\begin{align*}
(\star)&=E\bigg\{\sum_{k=\nu(j-1)+1}^{\nu(j)}\frac{|X_{k}|^{3}}{H_{k}^{3}}\bigg\vert
\Phi^{'''}\left(T_{k}(t)-\frac{\varepsilon_{k}X_{k}}{vH_{k}}\right)\bigg\vert\ind{A_{k}}\bigg\}\\
&\leq
\frac{c_{4}u_{n}}{m_{j}^{3}}\times\frac{v^{2}}{n}\times\left(\beta_{n-1}(4u,v)+\left(1-\frac{j-1}{n}\right)^{1/2}+M_{j}\right).
\end{align*}
In the other hand, for any integer $\nu(j-1)<k\leq\nu(j)$
$$
A_{k}^{c}\subset
B_{j}\triangleq\left\{\max_{\nu(j-1)<i\leq\nu(j)}\frac{\vert
R_{i}\vert}{m_{j}}>\frac{\vert
t-Y_{\nu(j-1)+1}\vert}{2M_{j}}\right\}.
$$
Since the set $A_{k}$ is $\F_{k}\vee\F_{\nu(j-1)}$, we have
\begin{align*}
(\star\star)&=E\bigg\{\sum_{k=\nu(j-1)+1}^{\nu(j)}\frac{|X_{k}|^{3}}{H_{k}^{3}}\bigg\vert
\Phi^{'''}\left(T_{k}(t)-\frac{\varepsilon_{k}X_{k}}{vH_{k}}\right)\bigg\vert\ind{A_{k}^{c}}\bigg\}\\
&\leq\|\Phi^{'''}\|_{\infty}E\left\{\sum_{k=\nu(j-1)+1}^{\nu(j)}\frac{|X_{k}|^{3}}{H_{k}^{3}}\ind{A_{k}^{c}}\right\}\\
&\leq u_{n}E\left\{\sum_{k=\nu(j-1)+1}^{\nu(j)}\frac{\sigma_{k}^{2}(X)}{H_{k}^{3}}\ind{A_{k}^{c}}\right\}\\
&\leq
u_{n}E\left\{\sum_{k=\nu(j-1)+1}^{\nu(j)}\frac{\sigma_{k}^{2}(X)}{H_{k}^{3}}\ind{B_{j}}\right\}.
\end{align*}
Since
\begin{equation}\label{inequality}
\sum_{k=\nu(j-1)+1}^{\nu(j)}\sigma_{k}^{2}(X)\leq\frac{2v^{2}}{n}\quad\textrm{a.s.}
\end{equation}
we derive
\begin{align*}
(\star\star)&\leq\frac{2u_{n}}{m_{j}^{3}}\times\frac{v^{2}}{n}\times\mu(B_{j})\\
&\leq\frac{2u_{n}}{m_{j}^{3}}\times\frac{v^{2}}{n}\times
\mu\left(\ds{\max_{\nu(j-1)<i\leq\nu(j)}}\vert R_{i}\vert>\frac{m_{j}\vert t-Y_{\nu(j-1)+1}\vert}{2M_{j}}\right)\\
&\leq\frac{2u_{n}}{m_{j}^{3}}\times\frac{v^{2}}{n}\times
E\left(\min\left\{1,\frac{4M_{j}^{2}}{m_{j}^{2}\vert
t-Y_{\nu(j-1)+1}\vert^{2}}
E\left(\max_{\nu(j-1)<i\leq\nu(j)}\vert R_{i}\vert^{2}\vert\F_{\nu(j-1)}\right)\right\}\right)\\
&\leq\frac{2u_{n}}{m_{j}^{3}}\times\frac{v^{2}}{n}\times
E\left(\min\left\{1,\frac{8M_{j}^{2}}{m_{j}^{2}\vert
t-Y_{\nu(j-1)+1}\vert^{2}}
E\left(\vert R_{\nu(j)}\vert^{2}\vert\F_{\nu(j-1)}\right)\right\}\right)\\
&\leq\frac{2u_{n}}{m_{j}^{3}}\times\frac{v^{2}}{n}\times
E\left(\min\left\{1,\frac{16M_{j}^{2}}{nm_{j}^{2}\vert
t-Y_{\nu(j-1)+1}\vert^{2}}\right\}\right)\quad\textrm{(using ($\ref{inequality}$))}\\
&\leq\frac{2u_{n}}{m_{j}^{3}}\times\frac{v^{2}}{n}\times\left(\beta_{n-1}(4u,v)+\left(1-\frac{j-1}{n}\right)^{1/2}+M_{j}\right)
\quad\textrm{(using Lemma $\ref{lemma-totalvariation}$)}.
\end{align*}
Thus there exists a positive constant $c_{5}$ such that
$$
(\star)+(\star\star)\leq\frac{c_{5}u_{n}}{m_{j}^{3}}\times\frac{v^{2}}{n}\times\left(\beta_{n-1}(4u,v)+\left(1-\frac{j-1}{n}\right)^{1/2}+M_{j}\right).
$$
Finally, we obtain the following estimate
\begin{align*}
S_{1}&\triangleq\sum_{k=1}^{n}E\left\{\frac{\vert
X_{k}\vert^{3}}{H_{k}^{3}}\bigg\vert\Phi^{'''}\left(T_{k}(t)-\frac{\varepsilon_{k}\theta_{k}}{vH_{k}}\right)\bigg\vert\right\}\\
&\leq c_{5}u_{n}\times\frac{v^{2}}{n}\times
\left(\beta_{n-1}(4u,v)\sum_{j=1}^{n}\frac{1}{m_{j}^{3}}
+\sum_{j=1}^{n}\frac{1}{m_{j}^{3}}\left(1-\frac{j-1}{n}\right)^{1/2}
+\sum_{j=1}^{n}\frac{M_{j}}{m_{j}^{3}}\right)\\
&\leq c_{5}u_{n}\times\frac{v^{2}}{n}\times
\left(\beta_{n-1}(4u,v)\frac{nv}{\sqrt{\theta^{2}-2u_{n}^{2}}}
+n\log n\right).
\end{align*}
Note that to obtain the above estimates of $S_{1}$, we have only
use the fact that the martingale difference sequence $X$ belongs
to the class $\L_{n}(u,v)$. Since the sequence $\sigma Z\triangleq
(\sigma_{1}(X)Z_{1},...,\sigma_{n}(X)Z_{n})$ belongs to
$\L_{n}(4u/\sqrt{2\pi},v)$ (with respect to the filtration
$\widetilde{\F}_{k}\triangleq\sigma(X_{1},...,X_{k},Z_{1},...,Z_{k})$),
we are able to reach a similar estimate for $S_{2}$:
$$
S_{2}\leq c_{6}u_{n}\times\frac{v^{2}}{n}\times
\left(\beta_{n-1}(16u/\sqrt{2\pi},v)\frac{nv}{\sqrt{\theta^{2}-2u_{n}^{2}}}
+n\log n\right)
$$
where $c_{6}$ is a positive constant. Using
($\ref{inequalitybeginning}$) and ($\ref{inequalityinduction}$),
there exist a positive constant $c$ such that
$$
\beta_{n}(u,v)\leq
c\,u_{n}\left(\frac{\beta_{n-1}(16u/\sqrt{2\pi},v)}{\sqrt{\theta^{2}-2u_{n}^{2}}}
+\frac{\log n}{v}\right)+\frac{3\theta}{v}.
$$
Putting
$$
D_{n}\triangleq\sup\left\{\frac{v\beta_{n}(u,v)}{u_{n}\log
n}\,\,;\,\,u\in\R_{+}^{\N^{\ast}},\,v>0\right\}
$$
and $\theta^{2}\triangleq \left(2+4c^{2}\right)u_{n}^{2}$, we
derive
$$
D_{n}\leq\frac{D_{n-1}}{2}+C
$$
where $C$ is a positive constant which does not depend on $n$.
Finally, we conclude that
$$
\limsup_{n\to+\infty}D_{n}\leq 2C.
$$
The proof of Theorem $\ref{maintheorem}$ is complete.\\
\\
{\em Proof of Theorem $\ref{cormaintheorem}$}. Let
$X=(X_{1},...,X_{n})$ in $\M_{n}(u)$. Following an idea by
Bolthausen, we are going to define a new martingale difference
sequence $\hat{X}$ which satisfies $V_{n}^{2}(\hat{X})=1$ a.s.
Denote $d_{1}=\|v_{n}^{2}V_{n}^{2}(X)-v_{n}^{2}\|_{1}$ and
$d_{\infty}=\|v_{n}^{2}V_{n}^{2}(X)-v_{n}^{2}\|_{\infty}$. The
letter $d$ will denote either $d_{1}$ or $d_{\infty}$. Consider
the random variables $X_{n+1},...,X_{n+[2d/u_{n}^{2}]+1}$ defined
as follows: Let $k=[(v_{n}^{2}+d-v_{n}^{2}V_{n}^{2})/u_{n}^{2}]$,
conditioned on $\F_{n+j-1}$, we assume
$$
X_{n+j}=\left\{\begin{array}{lll}
               \pm\,\,u_{n}\quad\textrm{w.p. $1/2$}&\textrm{for}&j\leq k\\
               \pm\,\,
               (v_{n}^{2}+d-v_{n}^{2}V_{n}^{2}-ku_{n}^{2})^{1/2}\quad\textrm{w.p. $1/2$}&\textrm{for}&j=k+1\\
               0&\textrm{else}&
               \end{array}\right.
$$
where $[\,.\,]$ denotes the integer part function and {\em w.p.}
is the abbreviation of {\em with probability}. In the sequel,
$\hat{n}$, $\hat{v}^{2}$, $\hat{V}^{2}$ and $\hat{S}$ denote
respectively $n+[2d/u_{n}^{2}]$, $v_{\hat{n}}^{2}(\hat{X})$,
$V_{\hat{n}}^{2}(\hat{X})$ and $S_{\hat{n}}(\hat{X})$. One can
easily check that $\hat{X}\triangleq (X_{1},...,X_{\hat{n}+1})$
belongs to $\M_{\hat{n}}(u)$ and $\hat{V}^{2}=1$ a.s. We have
$$
\Delta_{n}(X)\leq\sup_{t\in\R}\vert\mu(S_{n}/\hat{v}\leq
t)-\Phi(t)\vert+\sup_{t\in\R}\vert\Phi\left(\frac{v_{n}t}{\hat{v}}\right)-\Phi(t)\vert.
$$
Noting that $\hat{v}^{2}-v_{n}^{2}=d$ and using Lemma
$\ref{lemma-inequalities}$ with $k=2$ and $r=1$, if $d\triangleq
d_{1}$ there exist a positive constant $c$ such that
\begin{align*}
\Delta_{n}(X)&\leq 2\sup_{t\in\R}\vert\mu(\hat{S}/\hat{v}\leq
t)-\Phi(t)\vert+\frac{cd_{1}^{1/3}}{\hat{v}^{2/3}}+c\left(\frac{\hat{v}-v_{n}}{\hat{v}}\right)\\
&\leq 2\sup_{t\in\R}\vert\mu(\hat{S}/\hat{v}\leq
t)-\Phi(t)\vert+\frac{2cd_{1}^{1/3}}{v_{n}^{2/3}}.
\end{align*}
Similarly if $d\triangleq d_{\infty}$ then
\begin{align*}
\Delta_{n}(X)&\leq 2\sup_{t\in\R}\vert\mu(\hat{S}/\hat{v}\leq
t)-\Phi(t)\vert+\frac{cd_{\infty}^{1/2}}{\hat{v}}+c\left(\frac{\hat{v}-v_{n}}{\hat{v}}\right)\\
&\leq 2\sup_{t\in\R}\vert\mu(\hat{S}/\hat{v}\leq
t)-\Phi(t)\vert+\frac{2cd_{\infty}^{1/2}}{v_{n}}.
\end{align*}
Finally, applying Theorem $\ref{maintheorem}$ we derive
\begin{align*}
\Delta_{n}(X)&\leq 2\sup_{t\in\R}\vert\mu(\hat{S}/\hat{v}\leq
t)-\Phi(t)\vert+2c\min\left\{\frac{d_{1}^{1/3}}{v_{n}^{2/3}},\frac{d_{\infty}^{1/2}}{v_{n}}\right\}\\
&\leq 2L(u)\left(\frac{u_{\hat{n}}\log
\hat{n}}{\hat{v}}\right)+2c\min\left\{\frac{d_{1}^{1/3}}{v_{n}^{2/3}},\frac{d_{\infty}^{1/2}}{v_{n}}\right\}\\
&\leq 2L(u)\left(\frac{u_{n}\log
n}{v}\right)+2c\min\left\{\frac{d_{1}^{1/3}}{v_{n}^{2/3}},\frac{d_{\infty}^{1/2}}{v_{n}}\right\}\quad\textrm{if
$n$ is sufficiently large}.
\end{align*}
The proof of Theorem $\ref{cormaintheorem}$ is complete.\\
\\
{\em Proof of Lemma $\ref{lemma-inequalities}$}. Let $k>0$ and
$r\geq 1$, denote $\beta=\|E\left(\vert Y\vert^{k}\vert
X\right)\|_{r}$ and consider $q\in\R\cup\{\infty\}$ such that
$1/r+1/q=1$. Let $\lambda>0$ and $t$ be two real numbers we have
\begin{align*}
\mu\left(X+Y\leq t\right)&\geq\mu(X\leq t-\lambda,\,Y\leq t-X)\\
&=\mu(X\leq t-\lambda)-\mu(X\leq t-\lambda,\,Y\geq |t-X|)\\
&\geq\mu(X\leq t-\lambda)-E\left\{\ind{X\leq
t-\lambda}\,\mu(|Y|>|t-X\vert\,\vert X)\right\}.
\end{align*}
Since
\begin{align*}
E\left\{\ind{X\leq t-\lambda}\mu(|Y|>|t-X|\,\vert X)\right\}
&\leq E\left\{|t-X|^{-k}E(|Y|^{k}\vert X)\ind{X\leq t-\lambda}\right\}\\
&\leq\beta\|E\{\ind{X\leq t-\lambda}\vert t-X\vert^{-k}\}\|_{q}\\
&\leq\beta\lambda^{-k},
\end{align*}
we obtain
$$
\mu(X+Y\leq t)\geq\mu(X\leq
t-\lambda)-\beta\lambda^{-k}.
$$
Consequently
$$
\mu(X+Y\leq t)-\Phi(t)\geq\mu(X\leq
t-\lambda)-\Phi(t-\lambda)-\frac{\lambda}{\sqrt{2\pi}}-
\beta\lambda^{-k}
$$
and taking $\lambda=\left(\beta\sqrt{2\pi}\right)^{1/(k+1)}$,
there exists a positive constant $c$ such that
\begin{equation}\label{petite-inegalite1}
\delta(X+Y)\geq\delta(X)-c\beta^{1/(k+1)}.
\end{equation}
In the other hand
\begin{align*}
\mu(X+Y\leq t)&\leq\mu(X\leq t+\lambda)+\mu(X\geq t+\lambda, |Y|\geq |t-X|)\\
&=\mu(X\leq t+\lambda)+E\left\{\ind{X>t+\lambda}\,\mu(|Y|\geq
|t-X|\,\vert X)\right\}
\end{align*}
and
\begin{align*}
E\left\{\ind{X>t+\lambda}\,\mu(|Y|\leq |t-X|\,|X)\right\}
&\leq E\left\{\ind{X>t+\lambda}\,E(|Y|^{k}\,|X)\,|t-X|^{-k}\right\}\\
&\leq\beta\|E(\ind{X>t+\lambda}\,\vert t-X\vert^{-k})\|_{q}\\
&\leq\beta\lambda^{-k}.
\end{align*}
Consequently
$$
\mu(X+Y\leq t)\leq\mu(X\leq t+\lambda)+\beta\lambda^{-k}
$$
and
$$
\mu(X+Y\leq t)-\Phi(t)\leq\mu(X\leq
t+\lambda)-\Phi(t+\lambda)+\frac{\lambda}{\sqrt{2\pi}}+
\beta\lambda^{-k}.
$$
Taking $\lambda=(\beta\sqrt{2\pi})^{1/(k+1)}$, there exists a
positive constant $c^{'}$ such that
\begin{equation}\label{petite-inegalite2}
\delta(X+Y)\leq\delta(X)+c^{'}\beta^{1/(k+1)}.
\end{equation}
Combining ($\ref{petite-inegalite1}$) and
($\ref{petite-inegalite2}$) with Lemma 1 in Bolthausen
\cite{Bolth2} complete the proof of Lemma
$\ref{lemma-inequalities}$.\\
\\
{\em Proof of Theorem $\ref{rate-coboundary}$}. It suffice to
apply Inequality ($\ref{keyinequality1}$) of Lemma
$\ref{lemma-inequalities}$ with $k=p$, $r=1$ and $Y\triangleq
n^{-1/2}\left(g-g\circ T^{n}\right)$. The proof of Theorem $\ref{rate-coboundary}$ is complete.\\
\\
{\em Proof of Corollary $\ref{linear-process}$}. Since
$(\varepsilon_{j})_{j}$ is stationary, there exists a measure
preserving transformation $T$ such that
$\varepsilon_{j}=\varepsilon_{0}\circ T^{j}$. By Theorem 2 in
\cite{Volny93}, the condition ($\ref{condition-decomposition}$) is
necessary and sufficient to the existence of a function $g$ in
$L^{p}(\mu)$ such that
$$
X_{0}=m+g-g\circ T
$$
where $m\triangleq\varepsilon_{0}\sum_{k\in\Z}a_{k}$. Since $m$
satisfies the assumptions of Theorem $\ref{maintheorem}$, it
suffice to apply Theorem $\ref{rate-coboundary}$.
The proof of Corollary $\ref{linear-process}$ is complete.
\bibliographystyle{Paul1}
\bibliography{xbib}
\vspace{1cm}
Mohamed EL MACHKOURI, Lahcen OUCHTI\\
Laboratoire de Math\'ematiques Rapha\"el Salem\\
UMR 6085, Universit\'e de Rouen\\
Site Colbert\\
76821 Mont-Saint-Aignan, France\\
mohamed.elmachkouri@univ-rouen.fr\\
lahcen.ouchti@univ-rouen.fr
\end{document}